\documentclass[a4paper]{amsart}

\usepackage{algorithm}
\usepackage[noend]{algorithmic}
\usepackage{breakurl}
\usepackage[tableposition=below]{caption}
\usepackage{enumitem}
\usepackage{hyperref}
\usepackage{multirow}
\usepackage{tikz}
\usepackage{tikz-qtree}

\usetikzlibrary{positioning, calc}

\newtheorem{thm}{Theorem}[section]
\newtheorem{cor}[thm]{Corollary}
\newtheorem{rem}[thm]{Remark}
\theoremstyle{definition}
\newtheorem{defn}[thm]{Definition}

\newcommand{\N}{\mathbb{N}}
\newcommand{\Z}{\mathbb{Z}}
\newcommand{\Q}{\mathbb{Q}}
\newcommand{\F}{\mathbb{F}}

\DeclareMathOperator{\lm}{lm}
\DeclareMathOperator{\Lm}{Lm}
\DeclareMathOperator{\lc}{lc}
\DeclareMathOperator{\lt}{lt}
\DeclareMathOperator{\md}{mod}
\DeclareMathOperator{\Mon}{Mon}
\DeclareMathOperator{\NF}{NF}
\DeclareMathOperator{\tail}{tail}

\title{Gr\"obner Bases over Algebraic Number Fields}

\author[Dereje K. Boku]{Dereje Kifle Boku}
\address{Dereje Kifle Boku\\
Department of Mathematics\\
University of Kaiserslautern\\
Er\-win-Schr\"odinger-Str.\\
67663 Kaiserslautern\\
Germany}
\email{boku@mathematik.uni-kl.de}

\author{Wolfram Decker}
\address{Wolfram Decker\\
Department of Mathematics\\
University of Kaiserslautern\\
Erwin-Schr\"odinger-Str.\\
67663 Kaiserslautern\\
Germany}
\email{decker@mathematik.uni-kl.de}

\author{Claus Fieker}
\address{Claus Fieker\\
Department of Mathematics\\
University of Kaiserslautern\\
Erwin-Schr\"odinger-Str.\\
67663 Kaiserslautern\\
Germany}
\email{fieker@mathematik.uni-kl.de}

\author{Andreas Steenpass}
\address{Andreas Steenpass\\
Department of Mathematics\\
University of Kaiserslautern\\
Er\-win-Schr\"odinger-Str.\\
67663 Kaiserslautern\\
Germany}
\email{steenpass@mathematik.uni-kl.de}


\begin{document}

\begin{abstract}
Although Buchberger's algorithm, in theory, allows us to compute Gr\"obner
bases over any field, in practice, however, the computational efficiency
depends on the arithmetic of the ground field. Consider a field
$K = \Q(\alpha)$, a simple extension of $\Q$, where $\alpha$ is an algebraic
number, and let $f \in \Q[t]$ be the minimal polynomial of $\alpha$. In this
paper we present a new efficient method to compute Gr\"obner bases in
polynomial rings over the algebraic number field $K$. Starting from the ideas
of Noro \cite{B29}, we proceed by joining $f$ to the ideal to be considered,
adding $t$ as an extra variable. But instead of avoiding superfluous S-pair
reductions by inverting algebraic numbers, we achieve the same goal by applying
modular methods as in \cite{B5, BDFP, B45}, that is, by inferring information
in characteristic zero from information in characteristic $p > 0$. For suitable
primes $p$, the minimal polynomial $f$ is reducible over $\F_p$. This allows us
to apply modular methods once again, on a second level, with respect to the
factors of $f$. The algorithm thus resembles a divide and conquer strategy and
is in particular easily parallelizable. At current state, the algorithm is
probabilistic in the sense that, as for other modular Gr\"obner basis
computations, an effective final verification test is only known for
homogeneous ideals or for local monomial orderings. The presented timings show
that for most examples, our algorithm, which has been implemented in
\textsc{Singular} \cite{B25}, outperforms other known methods by far.
\end{abstract}

\maketitle

\section{Introduction}

From the theoretical point of view, Gr\"obner bases computations can be done
over any field by using Buchberger's algorithm (see, for example,
\cite{B4, B14, B6}). In particular, they can be performed over an algebraic
number field, but the computation is often inefficient if the arithmetic
operations in this field are used directly. Consider a simple extension
$K = \Q(\alpha)$ of $\Q$. Let $f \in \Q[t]$ be the minimal polynomial of
$\alpha$. The algebraic number field $K$ can be represented as the residue
class ring $\Q[t] / \langle f \rangle$, and a Gr\"obner basis computation over
$K$ can then be reduced to one over $\Q$ by joining $f$ to the ideal to be
considered. Unfortunately, this method is not satisfactory in view of
efficiency. One of the reasons for this is that over the field of rational
numbers, we often suffer from coefficient swell. Various methods to avoid this
have been investigated; the trace algorithm \cite{B26} and modular algorithms
\cite{B5, B45} are successful in this direction. But using these approaches, we
still have to deal with the complicated arithmetic in algebraic number fields,
in particular with the computation of inverses.

\vspace{0pt plus 2pt minus 2pt}
In this paper we present a new efficient method to compute Gr\"obner bases over
an algebraic number field. Starting from a polynomial ring over $\Q$ as
explained above, we apply the modular methods for computing Gr\"obner bases
discussed in \cite{B5, BDFP, B45} to pass to positive characteristic $p$.
Choosing a set $\mathcal{P}$ of suitable prime numbers, see
Definition~\ref{def:admissibleA}, the image $f_p$ of $f$ in $\F_p[t]$ is, for
$p \in \mathcal{P}$, reducible and square-free. We can thus again apply modular
methods, with respect to the factors $f_{1,p}, \ldots, f_{r_p,p}$ of $f_p$,
passing to the rings $\F_p[t] / \langle f_{i, p} \rangle$. As above, we avoid
computing in quotient rings by joining $f_{i, p}$ to the ideal to be
considered. Having computed the corresponding reduced Gr\"obner basis for each
of these factors, we first recombine the results to a set of polynomials $G_p$
over $\F_p[t] / \langle f_p \rangle$ using Chinese remaindering for
polynomials. In a second lifting step, the sets $G_p$, $p \in \mathcal{P}$, are
then used to reconstruct a set of polynomials $G$ over $\Q$, via Chinese
remaindering for integers and rational reconstruction. Finally, we test whether
$G$ is indeed the reduced Gr\"obner basis of the input ideal. If not, we
enlarge $\mathcal{P}$ and repeat the process.

\vspace{0pt plus 2pt minus 2pt}
In Section~\ref{sec:notation}, we introduce some notation which is used
throughout this article. The structure of the new method is outlined in
Section~\ref{sec:structure}. Since this method relies on the Chinese remainder
algorithm applied to different domains, we shortly recall the relevant
theoretical background in Section~\ref{sec:cra}. The core part of the proposed
algorithm is discussed in Section~\ref{sec:gbs}. Here we explain how modular
methods are applied on different levels and why our approach is considerably
faster than other known methods. The application of modular methods follows a
well-known scheme, see \cite{BDFP}. For reference, we recall the relevant parts
of this scheme in Section~\ref{sec:modular}. An illustrating example is given
in Section~\ref{sec:example}. Finally, Section~\ref{sec:timings} contains
remarks on the implementation of the new method in \textsc{Singular} \cite{B25}
and timings comparing it to other approaches. The benchmark problems which we
used for the timings are listed in the appendix.

\section{Notation}\label{sec:notation}

Let $K = \Q(\alpha)$ be an algebraic number field and let $f \in \Q[t]$ be the
minimal polynomial of the algebraic number $\alpha$. Then every element of $K$
can be written as a linear combination of elements in
$\lbrace 1, \alpha, \alpha^2, \ldots, \alpha^{d-1} \rbrace$ where $d = \deg f$.
Hence we may regard every element of $K$ as a polynomial in $\alpha$ with
coefficients in $\Q$. Let $X = \lbrace x_1, \ldots, x_n \rbrace$ be a set of
variables, and let $t$ be an extra variable. Consider the polynomial rings
$S = \Q(\alpha)[X]$, $T = \Q[X,t]$, and $\Q[t]$. Fix a global monomial ordering
$\succ_1$ on the monoid of monomials $\Mon(X)$ and consider the product
ordering $\succ_K := (\succ_1, \succ)$ on $\Mon(X, t)$, where $\succ$ is the
global ordering on $\Mon(t)$. Note that this implies $X^a \succ_K t^b$ for all
$a \in \N^n \setminus \lbrace (0, \ldots, 0) \rbrace$ and $b \in \N$.

Let $\widetilde{H} = \lbrace g_1(X,t), \ldots, g_s(X,t) \rbrace$ be a subset of
$T$, let $I \subseteq S$ be the ideal generated by
$H := \lbrace g_1(X,\alpha), \ldots, g_s(X,\alpha) \rbrace$, and let
$\widetilde{I} \subseteq T$ be the ideal generated by
$\widetilde{H} \cup \lbrace f \rbrace$. Furthermore, let
$\widetilde{G} \subseteq T$ be the reduced Gr\"obner basis (see
\cite[Definition~1.6.2]{B6}) of $\widetilde{I}$ w.r.t.\@ $\succ_K$. Let
$\varphi$ be the canonical homomorphism from $T$ to $S$ which leaves the $x_i$
fixed and maps $t$ to $\alpha$. We will show, in Theorem~\ref{thm:iso_redGB},
that the non-zero elements of $\varphi(\widetilde{G}) \subseteq S$ form the
reduced Gr\"obner basis of $I$ w.r.t.\@ $\succ_1$.

For a carefully chosen prime $p$ (see Definition~\ref{def:admissibleA}) which
does not divide any denominator of the coefficients of $f$ and
$g_1(X,t), \ldots, g_s(X,t)$, we consider the map from $\Q$ to $\F_p$ which
sends $\frac{a}{b}$ to $ab^{-1} \in \F_p$. Applying this map to the
coefficients, we write $f_p := (f \md p) \in \F_p[t]$ and
$\widetilde{I}_p := \langle g_1(X,t)_p, \ldots, g_s(X,t)_p, f_p \rangle
\subseteq \F_p[X,t]$.
Furthermore, for a polynomial $q \in S$ and a set $G \subseteq S$, we use the
following notation:

\noindent
$\lm(q)$: the \emph{leading monomial} of $q$, \\
$\Lm(G)$: the \emph{set of leading monomials} of the elements in $G$, \\
$\lc(q)$: the \emph{leading coefficient} of $q$, \\
$\lt(q)$: the \emph{leading term} or \emph{head} of $q$, \\
$\tail(q) := q-\lt(q)$: the \emph{tail} of $q$.

\section{Structure of the New Method}\label{sec:structure}

Noro \cite{B29} has presented a modified version of Buchberger's algorithm
which computes Gr\"obner bases over an algebraic number field using the
arithmetic in $\Q[t] / \langle f \rangle$. Instead of computing in the ring
$(\Q[t] / \langle f \rangle)[X]$, one might as well add the minimal polynomial
$f$ to the ideal to be considered and work over $\Q[X,t]$, see
Theorem~\ref{thm:iso_redGB}. In this situation, the elements of a reduced
Gr\"obner basis are, except $f$ itself, all monic in $(\Q[t])[X]$, that is,
they are of the form $X^a + (\text{lower terms})$, see the proof of
Theorem~\ref{thm:iso_redGB}. Noro noticed that during the execution of
Buchberger's algorithm, many (superfluous) intermediate basis elements of the
form $t^b X^a + (\text{lower terms})$ are computed before a monic element
$X^a + (\text{lower terms})$ is generated. Of course, each additional basis
element produces new S-pairs which usually make the subsequent computation
inefficient. Noro has resolved this problem by making each generated basis
element monic in $(\Q[t])[X]$ before it is added to the basis. For this, the
inverse of an algebraic number has to be computed which is in general
computationally expensive. Instead, we use a different approach to reduce the
number of basis elements which are computed before a monic element
$X^a + (\text{lower terms})$ is generated.

The new method computes the reduced Gr\"obner basis of the input ideal in three
steps: In the first step, for a suitable prime $p$ such that $f_p \in \F_p[t]$
is reducible and square-free, see Definition~\ref{def:admissibleA}, we compute
the reduced Gr\"obner basis $\widetilde{G}_p$ of $\widetilde{I}_p$ over $\F_p$
w.r.t.\@ $\succ_K$, as follows: Let $f_p = \prod_{1 \leq i \leq r_p} f_{i,p}$
be the irreducible factorization of $f_p$ over $\F_p$, with $r_p > 1$. Set
$\widetilde{I}_{i,p} := \langle \widetilde{H}_p \cup \lbrace f_{i,p} \rbrace
\rangle \subseteq \F_p[X,t]$.
For each $i \in \lbrace 1, \ldots, r_p \rbrace$, we compute the reduced
Gr\"obner basis $\widetilde{G}_{i,p}$ of $\widetilde{I}_{i,p}$. Using the
Chinese remainder algorithm for polynomials (see Algorithm~\ref{alg:cra}
below), we determine a set of polynomials
$\widetilde{G}_p \equiv \left( \widetilde{G}_{i,p} \setminus \lbrace f_{i,p}
\rbrace \right) \mod f_{i,p}$
which together with $f_p$ is the reduced Gr\"obner basis of $\widetilde{I}_p$
with high probability (see Remark~\ref{rem:high_probability}). Note that, at
this step of the algorithm, computing modulo the different factors of the
minimal polynomial $f_{i,p}$ (by adding them to the ideal
$\langle \widetilde{H}_p \rangle$) is, from the theoretical point of view, just
the same as computing modulo several prime numbers, see Section~\ref{sec:cra}.

In the second step, following \cite{B5, B45}, we use the Chinese remainder
algorithm for integers together with rational reconstruction to lift these
results to the reduced Gr\"obner basis $\widetilde{G}$ of $\widetilde{I}$. In
the last step, we lift $\widetilde{G}$ to a Gr\"obner basis $G$ of $I$ over $K$
by mapping $t$ to $\alpha$ (see Theorem~\ref{thm:iso_redGB}).

The idea of the algorithm is based on the concepts of modular methods and
univariate polynomial factorization over finite fields. For the former we need
the Chinese remainder theorem.

\section{Factorization and the Chinese Remainder Algorithm}\label{sec:cra}

The well-known Chinese remainder theorem is essential for our algorithm.

\begin{thm}[{\cite[Corollary~5.3]{B1}}]\label{thm:crt}
Let $R$ be a Euclidean domain and let $m_1, \ldots,\allowbreak m_r \in R$ be
coprime elements so that $\gcd(m_i, m_j) = 1$ for $0 \leq i < j \leq r$. Let
$m = m_1 \cdots m_r$ be the product of these elements. Then
$R / \langle m \rangle$ is isomorphic to the product ring
$R / \langle m_1 \rangle \times \ldots \times R / \langle m_r \rangle$ via the
isomorphism
\begin{align*}
R / \langle m \rangle
&\rightarrow R / \langle m_1 \rangle \times \ldots
\times R / \langle m_r \rangle \,, \\
a &\mapsto (a \md m_1, \ldots, a \md m_r) \,.
\end{align*}
\end{thm}

For our purpose, we need this theorem in the following two incarnations.

\begin{cor}
Let $p_1, \ldots, p_k$ be distinct prime numbers, and let $N = p_1 \cdots p_k$
be their product. Then we have the following isomorphism:
\[
\Z / \langle N \rangle \cong \F_{p_1} \times \ldots \times \F_{p_k} \,.
\]
\end{cor}

The second application of the Chinese remainder theorem refers to univariate
polynomial rings over finite fields.

\begin{cor}\label{cor:crt_poly}
Let $f_{1,p}, \ldots, f_{r_p,p} \in \F_p[t]$ be pairwise coprime polynomials,
and let $f_p = f_{1,p} \cdots f_{r_p,p}$ be their product. Then we have the
ring isomorphism
\[
\F_p[t] / \langle f_p \rangle
\cong \F_p[t] / \langle f_{1,p} \rangle \times \ldots
\times \F_p[t] / \langle f_{r_p,p} \rangle \,.
\]
\end{cor}

The proof of Theorem~\ref{thm:crt} is constructive (see
\cite[Theorem~5.2, Corollary~5.3]{B1}) and yields the Chinese remainder
algorithm. For reference, we state it here in the form of
Corollary~\ref{cor:crt_poly}, see Algorithm~\ref{alg:cra}.

\begin{algorithm}[h]
\caption{Chinese Remainder Algorithm (CRA) for polynomials}\label{alg:cra}

\begin{algorithmic}[1]
\REQUIRE $q_1, \ldots, q_{r_p} \in \F_p[t]$,
  $f_{1,p}, \ldots, f_{r_p,p} \in \F_p[t]$ pairwise coprime.
\ENSURE $g \in \F_p[t]$ such that $g \equiv q_i \mod f_{i,p}$ for
  $1\leq i\leq r_p$.

\STATE $g \longleftarrow 0$
\STATE $f_p \longleftarrow \prod_{1 \leq i \leq r_p} f_{i,p}$
\FOR {$i = 1, \ldots, r_p$}
  \STATE $h_i \longleftarrow \dfrac{f_p}{f_{i,p}}$
  \STATE by the Extended Euclidean Algorithm \cite[Algorithm~3.14]{B1}, compute
    $s_i, t_i \in \F_p[t]$ such that
    \[
    s_i h_i + t_i f_{i,p} = 1
    \]
  \STATE $c_i \longleftarrow \NF(q_i s_i, f_{i,p})$ \\
    ($c_i$ is the remainder in $\F_p[t]$ on dividing $q_i s_i$ by $f_{i,p}$)
  \STATE $g \longleftarrow g + c_i h_i$
\ENDFOR
\RETURN $g$
\end{algorithmic}
\end{algorithm}

\begin{rem}
\leavevmode
\begin{enumerate}[topsep=0pt,label=\alph*)]
\item
Since $c_i h_i \equiv 0 \mod f_{j,p}$ for $j \neq i$ and
$c_i h_i \equiv q_i s_i h_i \equiv q_i \mod f_{i,p}$, we have
\[
g \equiv c_i h_i \equiv q_i \mod f_{i,p} \,.
\]
Hence, the algorithm works correctly.

\item
Although stated here for $\F_p[t]$, Algorithm~\ref{alg:cra} works for
polynomial rings over any ground field.

\item
Instead of $q_1, \ldots, q_{r_p} \in \F_p[t]$, Algorithm~\ref{alg:cra} can also
be applied coefficient-wise to polynomials with coefficients in $\F_p[t]$.
\end{enumerate}
\end{rem}

\section{Gr\"obner Bases using Factorization and Modular Methods}%
\label{sec:gbs}

As Noro does (see \cite[Theorem~1]{B29}), we rely on the following result whose
proof we give for the lack of reference.

\begin{thm}\label{thm:iso_redGB}
Let $\widetilde{G}$ be the reduced Gr\"obner basis of $\widetilde{I}$ w.r.t.\@
$\succ_K$. Then $(\widetilde{G} \setminus \lbrace f \rbrace)\vert_{t=\alpha}$
is the reduced Gr\"obner basis of $I$ w.r.t.\@ $\succ_1$.
\end{thm}

Consider the ring homomorphism
\[
\varphi: T \longrightarrow S,\; t \longmapsto \alpha,\; x_i \longmapsto x_i \,.
\]
Since $\varphi$ is the identity map on $\Q[X]$, we get an isomorphism
\[
S \cong T / \langle f \rangle \,.
\]
Clearly, $\varphi(\widetilde{I}) = I$. We are now ready to prove
Theorem~\ref{thm:iso_redGB}.

\begin{proof}
Without loss of generality, we may assume that
$\widetilde{I} \neq \langle 1 \rangle$. Let
\[
\widetilde{G} = \lbrace m_1(X,t), \ldots, m_a(X,t), m_{a+1}(X,t) \rbrace
\]
be the reduced Gr\"obner basis of $\widetilde{I}$. We first prove that
$f \in \widetilde{G}$. Suppose $f \notin \widetilde{G}$. Then there exists a
non-zero non-constant polynomial $f' \in \widetilde{G} \cap \Q[t]$ with
$\deg f' < \deg f$. Hence
\[
I = \varphi(\widetilde{I})
= \langle \varphi(f'),
\varphi(\widetilde{G} \setminus \lbrace f' \rbrace) \rangle
= \langle 1 \rangle
\]
since $\varphi(f')$ is invertible in $S$. This implies
$\widetilde{I} = \langle 1 \rangle$, a contradiction. So, $f = m_i(X,t)$ for
some $i$, say $i = a+1$. Then we have
\begin{align*}
\varphi(\widetilde{G} \setminus \lbrace f \rbrace)
&= \lbrace m_1(X,\alpha), \ldots, m_a(X,\alpha) \rbrace \\
&= (\widetilde{G} \setminus \lbrace f \rbrace)\vert_{t=\alpha} =: G \,.
\end{align*}
The result follows easily once we show that the leading coefficient of
$m(X,t)$, considered as an element in the polynomial ring $\Q[t]$, is equal to
$1$ for all $m(X,t) \in \widetilde{G} \setminus \lbrace f \rbrace$. To prove
this statement, suppose there is an index $1 \leq j \leq a$ such that
$\lt(m_j(X,t)) = c \cdot X^\delta$ with $c \in \Q[t]$ and $\deg c > 0$.
Clearly, $c$ is monic. Write
\[
m_j(X,t) = c \cdot X^\delta + V(X,t)
\]
where $V(X,t) = \tail(m_j(X,t))$, which implies that $V(X,t)$ does not contain
any term divisible by $X^\delta$. We have $\deg c < \deg f$ and therefore
$\gcd(c, f) = 1$ since $f$ is irreducible. Thus, by the extended Euclidean
algorithm (see \cite[Algorithm~3.14]{B1}), there exist $a, b \in \Q[t]$ such
that $a \cdot c + b \cdot f = 1$. Considering the polynomial
$a \cdot m_j(X,t) + b \cdot f \cdot X^\delta$, we have
\begin{align*}
\langle \widetilde{G} \rangle
&\ni a \cdot m_j(X, t) + b \cdot f \cdot X^\delta \\
&= (a \cdot c + b \cdot f) \cdot X^\delta + a \cdot V(X,t) \\
&= X^\delta + a \cdot V(X,t) =: F(X,t) \,.
\end{align*}
But $\lt(F(X,t)) = X^\delta$ divides $c \cdot X^\delta = \lt(m_j(X,t))$ which
is a contradiction to the choice of $\widetilde{G}$.
\end{proof}

The notion of primes which are \emph{admissible of type~$\mathrm{A}$} w.r.t.\@
a monic irreducible polynomial, which is essential for our algorithm, is
defined as follows:
\begin{defn}\label{def:admissibleA}
Let $f \in \Q[t]$ be as given above. Let $p$ be a prime not dividing any
numerator or any denominator of the coefficients occurring in $f$. We say that
$p$ is \emph{admissible of type~$\mathrm{A}$} w.r.t.\@ $f$ if $f_p$ is
reducible and square-free over $\F_p$. In this case, we write $f_p$ as
$f_p = \prod_{1 \leq i \leq r_p} f_{i,p}$.
\end{defn}

For a non-zero polynomial $g \in T$ considered as a polynomial in $X$ over
$\Q[t]$, that is, $g \in (\Q[t])[X]$, let $S_g$ be the set of all distinct
coefficients (in $\Q[t]$) of $g$ of degree greater than or equal to 1. That is,
\[
S_g = \left\lbrace \lc_{\Q[t]}(u) \mid u \text{ is a term of } g \text{ with }
\deg(\lc_{\Q[t]}(u)) \geq 1 \right\rbrace \,.
\]

With notation as above, the notion of primes which are \emph{admissible of
type~$\mathrm{B}$} w.r.t.\@ a monic irreducible polynomial and a set of
polynomials is defined as follows:
\begin{defn}[Weak version]\label{def:admissibleB_weak}
Let $\widetilde{H} = \lbrace g_1(X,t),\allowbreak \ldots, g_s(X,t) \rbrace$ be
as given above. Let $p$ be a prime not dividing any numerator or any
denominator of the coefficients occurring in $\widetilde{H}$. We say that $p$
is \emph{admissible of type~$\mathrm{B}$} w.r.t.\@ $f$ and $\widetilde{H}$ if
$p$ is admissible of type~$\mathrm{A}$ w.r.t.\@ $f$ and if, for each $g$ in
$\widetilde{H}$, none of the elements in $S_g$ is divisible by any of the
factors of $f_p$ over $\F_p$.
\end{defn}

To see the relevance of this definition, consider the ideal
\[
J = \langle x^2+xy+t, x+y+t-1 \rangle =: \langle h_1, h_2 \rangle
\subseteq \Q[x,y,t]
\]
and the minimal polynomial $f = t^3+t+1$. If $p = 3$, then
$f_p \equiv (t-1)(t^2+t-1) =: f_{1,p} \cdot f_{2,p} \mod p$ and, using the
degree reverse lexicographic ordering with $x \succ y$, the reduced Gr\"obner
bases of the ideals $J_p + \langle f_{1,p} \rangle$ and
$J_p + \langle f_{2,p} \rangle$ in $\F_p[x,y,t]$ are $\lbrace 1 \rbrace$ and
$\lbrace t^2+t-1, y+1,x+t+1 \rbrace$, respectively. In this case,
Algorithm~\ref{alg:cra} cannot be applied since the sizes of these sets do not
fit. The calculation suggests that the reason for this is that the element
$t-1 \in S_{h_2}$ vanishes when reduced w.r.t.\@ the set
$\lbrace t-1, t^2+t-1 \rbrace$.

Next, consider the ideal
$J' = \langle x^2+xy+t, t^2x+y \rangle =: \langle g_1, g_2 \rangle$. Here, the
reduced Gr\"obner bases of the ideals $J'_p + \langle f_{1,p} \rangle$ and
$J'_p + \langle f_{2,p} \rangle$ are $\lbrace 1 \rbrace$ and
$\lbrace t^2+t-1, x+yt-t, y^2-1 \rbrace$, respectively. Again the sizes of
these sets do not coincide, hence, we still cannot apply
Algorithm~\ref{alg:cra}. Moreover, none of the coefficients in $S_{g_1}$ and
$S_{g_2}$ is divisible by either $f_{1,p}$ or $f_{2,p}$ which shows that the
condition in Definition~\ref{def:admissibleB_weak} is not sufficient. Indeed,
the element $t^2 \in S_{g_2}$ vanishes when reduced w.r.t.\@ the set
$\lbrace t^2+t-1, t-1 \rbrace$. Therefore, we may impose a stronger condition
by saying that for all $g \in \widetilde{H}$ none of the elements in $S_g$
vanishes when reduced w.r.t.\@ the set
$\lbrace f_{1,p}, \ldots, f_{r_p,p} \rbrace$ (in some order) and thus reduce
the probability that the reconstruction fails. In the following example we see
that this condition is still not sufficient.

Consider the ideal
$J'' = \langle x^2+xy+t, tx+y+t \rangle =:\langle k_1,k_2 \rangle$. The reduced
Gr\"obner bases of the ideals $J''_p + \langle f_{1,p} \rangle$ and
$J''_p + \langle f_{2,p} \rangle$ are $\lbrace t-1, x-1, y-1 \rbrace$ and
$\lbrace t^2+t-1, x+yt-y+t+1, y^2+yt+y+t-1 \rbrace$, respectively. Although
none of the elements in $S_{k_1}$ and $S_{k_2}$ vanishes when reduced w.r.t.\@
the set $\lbrace t^2+t-1, t-1 \rbrace$, and the sizes of these sets coincide,
we see that applying Algorithm~\ref{alg:cra} yields
$\lbrace t^2-t+1, x-1, y^2t^2+y^2t-y^2+yt^2+yt+t^2+t+1 \rbrace$ which is not
the desired result because the reduced Gr\"obner basis of
$J''_p + \langle f_p \rangle$ is $\lbrace t^2+t-1, y+1, x+t+1 \rbrace$. In
practice, however, it is very unlikely that this case happens. It is,
nevertheless, important to address this problem. A possible way to handle this
difficulty is to refine Definition~\ref{def:admissibleB_weak} as follows:

\begin{defn}[Strong version]\label{def:admissibleB_strong}
Let $f$ and $\widetilde{H} = \lbrace g_1(X,t), \ldots, g_s(X,t) \rbrace$ be as
given above. Let $p$ be an a prime which is admissible of type~$\mathrm{A}$
w.r.t.\@ $f$, and write $f = f_{1,p} \cdots f_{r_p,p}$ as in
Definition~\ref{def:admissibleA}. Suppose that $p$ does not divide any
numerator or any denominator of the coefficients occurring in $\widetilde{H}$.
For $i = 1, \ldots, r_p$, set
$\widetilde{I}_{i,p}
:= \langle \widetilde{H}_p \cup \lbrace f_{i,p} \rbrace \rangle$,
and let $\widetilde{G}_{i,p}$ be the reduced Gr\"obner basis of the ideal
$\widetilde{I}_{i,p}$. We say that $p$ is \emph{admissible of
type~$\mathrm{B}$} w.r.t.\@ $f$ and $\widetilde{H}$ if for all indices $i, j$
with $i \neq j$
\begin{enumerate}[label=\alph*)]
\item
the sizes of $\widetilde{G}_{i,p}$ and $\widetilde{G}_{j,p}$ coincide, and
\item
$\Lm(\widetilde{G}_{i,p} \setminus \lbrace f_{i,p} \rbrace)
= \Lm(\widetilde{G}_{j,p} \setminus \lbrace f_{j,p} \rbrace)$.
\end{enumerate}
\end{defn}

In the above examples, the prime number $3$ is not admissible of
type~$\mathrm{B}$ w.r.t.\@ $t^3+t+1$ and the generators of each of the ideals
$J$, $J'$ and $J''$ in the sense of Definition~\ref{def:admissibleB_strong}.
This is because in the first two cases, both conditions of this definition are
violated whereas in the third case, the second condition is not satisfied. For
the rest of our discussion we use the strong version of this definition.

We now turn our attention to the notion of \emph{lucky primes}:
\begin{defn}[{\cite{B45}}]
Let $\widetilde{I}$ be an ideal given as above and let $p$ be a prime number.
Furthermore, let $\widetilde{G}$ be the reduced Gr\"obner basis of
$\widetilde{I}$ and let $\widetilde{G}_{p}$ be the reduced Gr\"obner basis of
$\widetilde{I}_{p}$. Then $p$ is called \emph{lucky} for $\widetilde{I}$ if and
only if $\Lm(\widetilde{G}_p) = \Lm(\widetilde{G})$. Otherwise $p$ is called
\emph{unlucky} for $\widetilde{I}$.
\end{defn}

Since $f$ is independent of $X$, we get, by Corollary~\ref{cor:crt_poly}, the
isomorphism
\[
\F_p[X,t] / \langle f_p \rangle
\cong \F_p[X,t] / \langle f_{1,p} \rangle \times \ldots
\times \F_p[X,t] / \langle f_{r_p,p} \rangle \,.
\]

\begin{rem}\label{rem:high_probability}
Let $\widetilde{I}$, $\widetilde{H}$, and $f$ be as above. Let $p$ be a prime
which is both admissible of type~$\mathrm{B}$ w.r.t.\@ $f$ and $\widetilde{H}$
as well as lucky for $\widetilde{I}$. We work over $\F_p[X,t]$ equipped with
the product ordering $\succ_K$. Suppose a set of polynomials $\widetilde{G}_p$
is the reduced Gr\"obner basis of the ideal $\widetilde{I}_p$. For
$i =1, \ldots, r_p$, set
$S_i := (\widetilde{G}_p \setminus \lbrace f_p \rbrace) \md f_{i,p}
\subseteq \F_p[X,t] / \langle f_{i,p} \rangle$.
Then for each $i$, the set $S_i \cup \lbrace f_{i,p} \rbrace$ is the reduced
Gr\"obner basis of the ideal $\widetilde{I}_{i,p}$ (as in
Definition~\ref{def:admissibleB_strong}) with high probability. Conversely, let
$\widetilde{G}_{i,p}$ be the reduced Gr\"obner basis of $\widetilde{I}_{i,p}$.
Let $\widetilde{G}_p'$ be the set of polynomials that is obtained by applying
Algorithm~\ref{alg:cra} coefficient-wise to the input
\[
\left( (\widetilde{G}_{1,p} \setminus \lbrace f_{1,p} \rbrace, \ldots,
\widetilde{G}_{r_p,p} \setminus \lbrace f_{r_p,p} \rbrace),\,
(f_{1,p}, \ldots, f_{r_p,p}) \right) \,.
\]
Then the set $\widetilde{G}_p' \cup \lbrace f_{p} \rbrace$ is the reduced
Gr\"obner basis of the ideal $\widetilde{I}_p$ with high probability. Hence, we
have $\widetilde{G}_p' \cup \lbrace f_{p} \rbrace = \widetilde{G}_p$ with high
probability.
\end{rem}

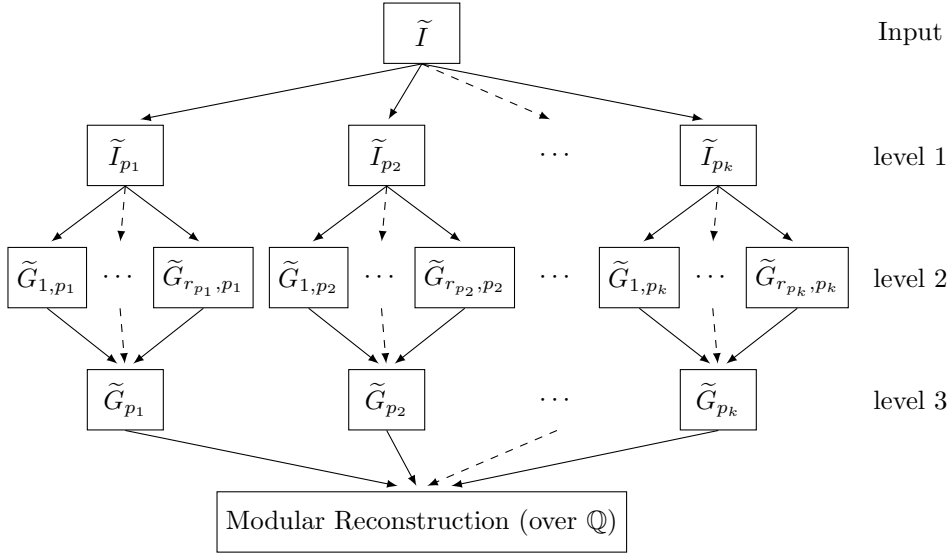
\begin{figure*}
\centering
\begin{tikzpicture}[
  -latex,
  every node/.style={
    minimum height=0.8cm,
    minimum width=1cm,
  },
  box/.style={
    draw,
    anchor=north,
  },
  nobox/.style={
    anchor=north,
  },
]

\newlength{\leveldist}
\setlength{\leveldist}{1.62cm}

\begin{scope}[
  level distance=\leveldist,
  sibling distance=0.2cm,
  edge from parent/.style={
    draw,
    edge from parent path={
      let \p1=($(\tikzchildnode.north)+(0pt,2pt)$),
      \p2=($(\tikzparentnode.south)!0.95!(\p1)$)
      in (\tikzparentnode.south) -- (\x2,\y1)
    }
  }
]
\Tree
[ .\node [box] (I) {$\widetilde{I}$};
  [ .\node [box] (I1) {$\widetilde{I}_{p_1}$};
      \node [box] (Gp11) {$\widetilde{G}_{1,p_1}$};
      \edge [dashed];
      \node [nobox, minimum width=0cm, inner sep=0pt] (cdots21) {$\cdots$};
      \node [box] (Gp1r) {$\widetilde{G}_{r_{p_1},p_1}$};
  ]
  [ .\node [box] (I2) {$\widetilde{I}_{p_2}$};
      \node [box] (Gp21) {$\widetilde{G}_{1,p_2}$};
      \edge [dashed];
      \node [nobox, minimum width=0cm, inner sep=0pt] (cdots22) {$\cdots$};
      \node [box] (Gp2r) {$\widetilde{G}_{r_{p_2},p_2}$};
  ]
  \edge [dashed];
  [ .\node [nobox, minimum width=0cm] (cdots1) {$\cdots$};
      \edge [draw=none];
      \node [nobox, minimum width=0cm, inner sep=0pt] (codts23) {$\cdots$};
  ]
  [ .\node [box] (Ik) {$\widetilde{I}_{p_k}$};
      \node [box] (Gpk1) {$\widetilde{G}_{1,p_k}$};
      \edge [dashed];
      \node [nobox, minimum width=0cm, inner sep=0pt] (cdots24) {$\cdots$};
      \node [box] (Gpkr) {$\widetilde{G}_{r_{p_k},p_k}$};
  ]
]
\end{scope}

\node [box, below=2\leveldist of I1.north] (G1) {$\widetilde{G}_{p_1}$};
\node [box, below=2\leveldist of I2.north] (G2) {$\widetilde{G}_{p_2}$};
\node [nobox, below=2\leveldist of cdots1.north, minimum width=0pt] (cdots3)
  {$\cdots$};
\node [box, below=2\leveldist of Ik.north] (Gk) {$\widetilde{G}_{p_k}$};
\node [box, below=4\leveldist of I.north] (bottom)
  {Modular Reconstruction (over $\Q$)};

\newcommand{\arrow}[3]{\draw [#3] let \p1=($(#2.north)+(0pt,2pt)$),
      \p2=($(#1.south)!0.90!(\p1)$)
      in (#1.south) -- (\x2,\y1);}
\arrow{Gp11}{G1}{};
\arrow{cdots21}{G1}{dashed};
\arrow{Gp1r}{G1}{};
\arrow{Gp21}{G2}{};
\arrow{cdots22}{G2}{dashed};
\arrow{Gp2r}{G2}{};
\arrow{Gpk1}{Gk}{};
\arrow{cdots24}{Gk}{dashed};
\arrow{Gpkr}{Gk}{};
\draw (G1.south) -- ($(bottom.north)+(-9pt,2pt)$);
\draw (G2.south) -- ($(bottom.north)+(-2pt,2pt)$);
\draw [dashed] (cdots3.south) -- ($(bottom.north)+(2pt,2pt)$);
\draw (Gk.south) -- ($(bottom.north)+(11pt,2pt)$);

\node [right=0.2cm of Gpkr] (level2) {level 2};
\node at (level2|-I) {Input};
\node at (level2|-Ik) {level 1};
\node at (level2|-Gk) {level 3};
\end{tikzpicture}
\caption{General scheme for the new algorithm}\label{fig:scheme}
\end{figure*}

The main innovation of our new algorithm, which is illustrated in
Figure~\ref{fig:scheme}, is as follows: Instead of computing the reduced
Gr\"obner bases at level~1, our algorithm computes them at level~2. For the
primes satisfying the conditions in Definition~\ref{def:admissibleB_strong}
(and only for those), the Chinese remainder algorithm for polynomials then
combines these results at level~3. The ideals
$\langle \widetilde{G}_{p_i} \rangle$ at this level are expected to be the same
as the ideals $\widetilde{I}_{p_i}$ at level~1 with high probability (see
Remark~\ref{rem:high_probability}). The remaining parts of the computation are
carried out in the same way as in the modular algorithms described in
\cite{B45}.

Now we give a brief description of the new algorithm. In the beginning,
randomly choose a set $\mathcal{P}$ of prime numbers which are admissible of
type~$\mathrm{A}$ w.r.t.\@ $f$. At level~2, given a prime $p \in \mathcal{P}$,
factorize $f \in \Q[t]$ over $\F_p$ and compute, for each $i$, the reduced
Gr\"obner basis $\widetilde{G}_{i,p}$ of the ideal $\widetilde{I}_{i,p}$
corresponding to the $i$-th factor. If the prime $p$ is admissible of
type~$\mathrm{B}$ w.r.t.\@ $f$ and $\widetilde{H}$, then lift these results via
Chinese remaindering for polynomials (at level~3) to obtain the reduced
Gr\"obner basis $\widetilde{G}_p$ of $\widetilde{I}_p$ with high probability.
Repeat this process for every prime $p \in \mathcal{P}$ which is admissible of
type~$\mathrm{B}$, in the same way as in the modular algorithms in \cite{B45}.

The main reason why the method to compute Gr\"obner bases over algebraic number
fields described above is faster than other known methods, see
Section~\ref{sec:timings}, is that factorizing the minimal polynomial $f$ in
positive characteristic allows us to compute in rings with minimal polynomials
of degree much less than $\deg f$: Experiments have shown that the performance
of Gr\"obner basis computations over simple algebraic extensions depends
heavily on the degree of the minimal polynomial. Additionally, the computations
are carried out over finite fields which avoids the problem known as
coefficient swell, and we do not directly use the computationally expensive
arithmetic in $K$. Finally, the new method is a priori easily parallelizable.

\section{Modular Algorithms}\label{sec:modular}

To compute the reduced Gr\"obner basis of the ideal $\widetilde{I}$, the
modular algorithm described in \cite{B45} first chooses a set of primes
$\mathcal{P}$ and computes the reduced Gr\"obner bases $\widetilde{G}_p$ of
$\widetilde{I}_p$ for each $p \in \mathcal{P}$. It then uses the Chinese
remainder algorithm and rational reconstruction to obtain the reduced Gr\"obner
basis $\widetilde{G}$ over $\Q$ with high probability. Finally, it verifies the
correctness of the result obtained in this way. One of the problems after
computing the set of reduced Gr\"obner bases
$\mathcal{GP} := \lbrace \widetilde{G}_p \mid p \in \mathcal{P} \rbrace$ is
that $\mathcal{P}$ may contain unlucky primes. To deal with such unlucky
primes, the following method is used, see \cite{BDFP}:

\textsc{DeleteUnluckyPrimesSB} (\cite{B45}):
\textit{We define an equivalence relation on
$(\mathcal{GP},\allowbreak \mathcal{P})$ by
\[
(\widetilde{G}_p, p) \sim (\widetilde{G}_q, q)
:\Longleftrightarrow \Lm(\widetilde{G}_p) = \Lm(\widetilde{G}_q) \,.
\]
Then the equivalence class of largest cardinality\footnote{Here, we have to use
a weighted cardinality count if Algorithm~\ref{alg:modStd_variant} requires
more than one round of the loop, see \cite[Remark~5.7]{BDFP}.} is stored in
$(\mathcal{GP}, \mathcal{P})$, the others are deleted}.

Now, all $\widetilde{G}_p$, $p \in \mathcal{P}$, have the same set of leading
monomials. Hence, we can apply the Chinese remainder algorithm for integers and
the rational reconstruction algorithm to the coefficients of the Gr\"obner
bases in $\mathcal{GP}$ to obtain a reduced Gr\"obner basis $\widetilde{G}$ of
$\widetilde{I}$ with high probability. Since we cannot check, however, whether
$\mathcal{P}$ is sufficiently large, a final verification step is needed. Since
this may be expensive, especially if
$\widetilde{I} \neq \langle \widetilde{G} \rangle$, we first perform a test in
positive characteristic:

\textsc{pTestSB} (\cite{B45}):
\textit{We randomly choose a prime $p \notin \mathcal{P}$ which is admissible
of type~$\mathrm{B}$ w.r.t.\@ $f$ and $\widetilde{H}$. We test if including
this prime in the set $\mathcal{P}$ would improve the result. That is,
explicitly test whether $\widetilde{I}$ reduces to zero w.r.t\@ $\widetilde{G}$
mapped to $\F_p[X, t]$, and vice-versa, whether $\widetilde{G}$
mapped to $\F_p[X, t]$ reduces to zero w.r.t.\@ $\widetilde{G}_p$.}

The advantage of this test is that it accelerates the algorithm enormously.
Algorithm~\ref{alg:modStd_variant} is a modified version of Algorithm~1 in
\cite{B45} (which is implemented in \textsc{Singular} \cite{B25} in the library
\texttt{modstd.lib} \cite{GLP}), in the sense that we do apply modular methods
not only once, but twice, where the second application is with respect to the
factors of the minimal polynomial $f$.

\begin{algorithm}[t]
\caption{Modified modular Gr\"obner bases algorithm over $\Q$}%
\label{alg:modStd_variant}

\begin{algorithmic}[1]
\REQUIRE an ideal
  $\widetilde{I} = \langle \widetilde{H}, f \rangle \subseteq T = \Q[X,t]$
  where $\widetilde{H} = \lbrace g_1(X,t), \ldots, g_s(X,t) \rbrace$ and
  $f \in \Q[t]$ is irreducible.
\ENSURE $\widetilde{G} \subseteq T$, a Gr\"obner basis of $\widetilde{I}$
  w.r.t.\@ $\succ_K$.

\STATE\label{line:modStd_variant_primes}
  choose $\mathcal{P}$, a set of random primes which are admissible of
  type~$\mathrm{A}$ w.r.t.\@ $f$
\STATE $\mathcal{GP} \longleftarrow \lbrace \rbrace$
\LOOP
  \FOR {$p \in \mathcal{P}$}\label{line:modStd_variant_for1}
    \STATE factorize $f_p \in \F_p[t]$ into irreducible factors
      $f_p = \prod_{1 \leq i \leq r_p} f_{i,p}$
    \FOR {$i = 1, \ldots, r_p$}\label{line:modStd_variant_for2}
      \STATE $\widetilde{I}_{i,p} \longleftarrow \langle \widetilde{H}_p \cup
        \lbrace f_{i,p} \rbrace \rangle \subseteq \F_p[X,t]$
      \STATE compute the reduced Gr\"obner basis $\widetilde{G}_{i,p}$ of
        $\widetilde{I}_{i,p}$ w.r.t.\@ $\succ_K$
    \ENDFOR
    \IF {$p$ is admissible of type~$\mathrm{B}$ w.r.t.\@ $f$ and
        $\widetilde{H}$ over $\F_p$}
      \STATE apply Algorithm~\ref{alg:cra} coefficient-wise to the input
        $\Bigl( \bigl( \widetilde{G}_{1,p} \setminus \lbrace f_{1,p} \rbrace,\,
        \ldots,\allowbreak\, \widetilde{G}_{r_p,p} \setminus \lbrace f_{r_p,p}
        \rbrace \bigr),\; \bigl(f_{1,p},\, \ldots,\allowbreak\, f_{r_p,p}
        \bigr) \Bigr)$
        to obtain a set of polynomials $\widetilde{G}_p \subseteq \F_p[X,t]$
      \STATE $\widetilde{G}_p \longleftarrow \widetilde{G}_p \cup \lbrace f_{p}
        \rbrace$
    \ELSE
      \STATE $\widetilde{G}_p \longleftarrow 0$
    \ENDIF
    \STATE $\mathcal{GP} \longleftarrow \mathcal{GP} \cup \lbrace
      \widetilde{G}_p \rbrace$
  \ENDFOR
  \STATE $(\mathcal{GP}, \mathcal{P}) \longleftarrow
    \textsc{DeleteUnluckyPrimesSB}(\mathcal{GP}, \mathcal{P})$
  \STATE\label{line:modStd_variant_farey}
    lift $(\mathcal{GP}, \mathcal{P})$ to $\widetilde{G} \subseteq T$ by
    applying the Chinese remainder algorithm and the Farey rational map
  \IF {$\textsc{pTestSB}(\widetilde{I}, \widetilde{G}, \mathcal{P})$}%
     \label{line:modStd_variant_finalTest_start}
    \IF {$\widetilde{I}$ reduces to zero w.r.t.\@ $\widetilde{G}$}%
       \label{line:modStd_variant_finalTest_end}
      \IF {$\widetilde{G}$ is the reduced Gr\"obner basis of
          $\langle \widetilde{G} \rangle$}\label{line:modStd_variant_finalTest}
        \RETURN $\widetilde{G}$
      \ENDIF
    \ENDIF
  \ENDIF
  \STATE enlarge $\mathcal{P}$
\ENDLOOP
\end{algorithmic}
\end{algorithm}

Now, taking Theorem~\ref{thm:iso_redGB} into account, we can compute a
Gr\"obner basis of an ideal in $K[X] = \Q(\alpha)[X]$ as in
Algorithm~\ref{alg:nfmodStd}: We first map $\alpha$ to $t$ and join the minimal
polynomial $f \in \Q[t]$ to the ideal to be considered. Then, after applying
Algorithm~\ref{alg:modStd_variant}, we only need to map $t$ back to $\alpha$ to
get a Gr\"obner basis of the input ideal.

\begin{algorithm}
\caption{Modular Gr\"obner basis algorithm over $K = \Q(\alpha)$
(\texttt{nfmodStd})}%
\label{alg:nfmodStd}

\begin{algorithmic}[1]
\REQUIRE $I = \langle g_1(X,\alpha), \ldots, g_s(X,\alpha) \rangle
  \subseteq S = K[X]$.
\ENSURE $G \subseteq S$, a Gr\"obner basis of $I$ w.r.t.\@ $\succ_1$.

\STATE map $I$ to $\langle \widetilde{H} \rangle$ via the map sending $\alpha$
  to $t$
\STATE $\widetilde{I} \longleftarrow \langle \widetilde{H} \rangle
  + \langle f \rangle$
\STATE call Algorithm~\ref{alg:modStd_variant} to compute the reduced Gr\"obner
  basis $\widetilde{G}$ of $\widetilde{I}$ w.r.t.\@
  $\succ_K =(\succ_1,\allowbreak \succ)$
\STATE lift $\widetilde{G}$ to $G$ via the map sending $t$ to $\alpha$
\RETURN $G$
\end{algorithmic}
\end{algorithm}

Algorithm~\ref{alg:modStd_variant} is probabilistic in the sense that the test
in lines \ref{line:modStd_variant_finalTest_start} to
\ref{line:modStd_variant_finalTest_end} does not guarantee that
$\langle \widetilde{G} \rangle = \widetilde{I}$. If $I$ is homogeneous,
however, the result $G$ of Algorithm~\ref{alg:nfmodStd} can be verified along
the lines of \cite[Theorem~7.1]{B5}. With this test included,
Algorithm~\ref{alg:nfmodStd} is deterministic.

\begin{rem}
Some parts of Algorithm~\ref{alg:modStd_variant} are inherently parallelizable.
In the current implementation, see Section~\ref{sec:timings}, we could easily
take advantage of this thanks to \textsc{Singular}'s parallel framework. We
have, first of all, parallelized the for-loop starting in
line~\ref{line:modStd_variant_for1}. This corresponds to the modular
computations on level~1, see Figure~\ref{fig:scheme}. Besides this, we also
make use of parallelization for the selection of primes in
line~\ref{line:modStd_variant_primes}, for the application of the Farey
rational map in line~\ref{line:modStd_variant_farey}, and for the final test in
line~\ref{line:modStd_variant_finalTest}. The for-loop starting in
line~\ref{line:modStd_variant_for2}, which corresponds to the modular
computations on level~2, is inherently parallelizable as well, but experiments
have shown that a parallel implementation of this step does not yield any
further speedup for our test cases.
\end{rem}

\section{Example}\label{sec:example}

The following example illustrates how the new algorithm works:

Consider the ideal
$I = \langle x^2 + ay, axy - x + a \rangle \subset \Q(a)[x,y]$ where $a$ is a
zero of the polynomial $f = t^2 + 1 \in \Q[t]$. A \textsc{Singular} computation
shows that the reduced Gr\"obner basis of $I$ with respect to the degree
reverse lexicographical ordering (\verb+dp+ in \textsc{Singular}) with
$x \succ y$ is
\[
\lbrace y^2+ax+ay,\; xy+ax+1,\; x^2+ay \rbrace \,.
\]

In the following, we show how this basis is obtained using our method: At
level~1, let us choose $k = 2$ with $p_1 = 5$ and $p_2 = 13$. At level~2, we
have $f_{p_1} \equiv (t-2)(t+2) \mod p_1$ and
$f_{p_2} \equiv (t-5)(t+5) \mod p_2$. Now, corresponding to each factor, we
compute, using {\sc{Singular}}, the reduced Gr\"obner bases of the following
ideals:
\begin{align*}
\widetilde{I}_{1,p_1} = \langle x^2 + ty , txy -x+t, t-2 \rangle \,, \\
\widetilde{I}_{2,p_1} = \langle x^2 + ty , txy -x+t, t+2 \rangle \,, \\
\widetilde{I}_{1,p_2} = \langle x^2 + ty , txy -x+t, t-5 \rangle \,, \\
\widetilde{I}_{2,p_2} = \langle x^2 + ty , txy -x+t, t+5 \rangle \,.
\end{align*}
\begin{verbatim}
> ring r = 5, (x,y,t), (dp(2),dp(1));
> ideal I1p1 = x2+ty, txy-x+t, t-2;
> ideal I2p1 = x2+ty, txy-x+t, t+2;
> option(redSB);
> ideal S1 = std(I1p1);
> S1;
  S1[1]=t-2
  S1[2]=y2+2x+2y
  S1[3]=xy+2x+1
  S1[4]=x2+2y
> ideal S2 = std(I2p1);
> S2;
  S2[1]=t+2
  S2[2]=y2-2x-2y
  S2[3]=xy-2x+1
  S2[4]=x2-2y
\end{verbatim}

The Chinese remainder algorithm for polynomials combines these results at
level~3 to obtain the reduced Gr\"obner basis of $\widetilde{I}_{p_1}$ with
high probability, as follows:
\begin{verbatim}
> list l = S1, S2;
> list m = t-2, t+2;
  // CRA for polynomials (coefficient-wise):
> ideal G1p1 = chinrempoly(l, m);
> Gp1;
  Gp1[1]=t2+1
  Gp1[2]=y2+xt+yt
  Gp1[3]=xy+xt+1
  Gp1[4]=x2+yt
\end{verbatim}

Similarly, the reduced Gr\"obner basis of $\widetilde{I}_{p_2}$, with high
probability, is
\begin{verbatim}
> Gp2;
  Gp2[1]=t2+1
  Gp2[2]=y2+xt+yt
  Gp2[3]=xy+xt+1
  Gp2[4]=x2+yt
\end{verbatim}

It is not hard to see that the primes $p_1$ and $p_2$ are admissible of
type~$\mathrm{B}$ w.r.t.\@ $f$ and
$\widetilde{H} = \lbrace x^2 + ty , txy - x + t \rbrace$. Furthermore, it is
also clear that they are lucky primes for
$\widetilde{I} = \langle\widetilde{H}, f \rangle$. At this point we have to
change the current base ring in \textsc{Singular} to characteristic zero in
order to apply the Chinese remainder algorithm for integers and to pull the
modular coefficients back to the rational numbers.
\begin{verbatim}
/* Chinese remaindering for integers */
> ring s = 0, (x,y,t), (dp(2),dp(1));
> list l = imap(r, Gp1), imap(r, Gp2);
> intvec m = 5, 13;
> ideal j = chinrem(l, m);
> j;
  j[1]=t2+1
  j[2]=y2+xt+yt
  j[3]=xy+xt+1
  j[4]=x2+yt
/* rational reconstruction */
> j = farey(j, 5*13);
> j;
  j[1]=t2+1
  j[2]=y2+xt+yt
  j[3]=xy+xt+1
  j[4]=x2+yt
\end{verbatim}

Note that the computed result already coincides with the reduced Gr\"obner
basis stated above. To simplify the presentation, we therefore skip some of the
steps in Algorithm~\ref{alg:modStd_variant}, such as the final test. However,
we have to map the result back to the ring $\Q(a)[x,y]$ in \textsc{Singular}:
\begin{verbatim}
> ring sr = (0,a), (x,y,t), (dp(2),dp(1));
> minpoly = a2+1;
> ideal G = imap(s, j);
> G = subst(G, t, a);
> G = simplify(G, 2); // erase the zero entries
> G; // G is the reduced Groebner basis of I
  G[1]=y2+ax+ay
  G[2]=xy+ax+1
  G[3]=x2+ay
\end{verbatim}
Thus we get the same result as the one we mentioned at the beginning.

\section{Implementation and Timings}\label{sec:timings}

We implemented Algorithm~\ref{alg:nfmodStd} in \textsc{Singular} in the library
\texttt{nfmodstd.lib} \cite{nfmodstd} and compared its performance against the
implementation of \cite[Algorithm~1]{B45} in the \textsc{Singular} library
\texttt{modstd.lib} (the command is \texttt{modStd}), the \textsc{Singular}
command \texttt{std}, and the Magma \cite{magma, B42} command
\texttt{GroebnerBasis}. For \texttt{modStd}, we added the minimal polynomial
$f$ to the given input ideal $I$ (considered as an ideal in a polynomial ring
over a polynomial ring) and computed the reduced Gr\"obner basis of the ideal
$\widetilde{I} = \langle \widetilde{H} \rangle + \langle f \rangle$ w.r.t.\@
$\succ_K$. For \texttt{GroebnerBasis} and \texttt{std}, we computed the reduced
Gr\"obner basis of the ideal $I$ over an algebraic number field with the
minimal polynomial $f$. Note that the implementation of our algorithm is
internally linked with the existing implementation of Algorithm~1 in
\cite{B45}.

We have nine benchmark problems to demonstrate the superiority of our new
algorithm (see appendix). The cyclic ideal $C_n$ in $n$ variables has become a
benchmark problem for Gr\"obner basis techniques. For our algorithm, we have
replaced the coefficients of this ideal by a random element in $\Q(a)$ where
$a$ is an algebraic number (see, for example, the ideal \texttt{I6} in the
appendix). Some of the benchmark problems are chosen from \cite{B5, B29} (the
ideals \texttt{I1} and \texttt{I2} are from \cite{B5}, \texttt{I6} and
\texttt{I7} are from \cite{B29}) where the coefficients are replaced by a
random algebraic number. The minimal polynomials, selected for our
computations, are:
{\allowdisplaybreaks[1]
\begin{align*}
m_1 &= a^2+1 \,, \\
m_2 &= a^5+a^2+2 \,, \\
m_3 &= a^7-7a+3 \,, \\
m_4 &= a^6+a^5+a^4+a^3+a^2+a+1 \,, \\
m_5 &= a^{12}-5a^{11}+24a^{10}-115a^9+551a^8-2640a^7 \\*
    &+12649a^6-2640a^5+551a^4-115a^3+24a^2-5a+1 \,, \\
m_6 &= a^2+5a+1 \,, \\
m_7 &= a^8-16a^7+19a^6-a^5-5a^4+13a^3-9a^2+13a \\*
    &+17 \,, \text{ and} \\
m_8 &= a^7+10a^5+5a^3+10a+1 \,.
\end{align*}}

With respect to these minimal polynomials, timings are conducted by using
\textsc{Singular}~4.0.2 and Magma \hbox{V2.21-2} on a Dell PowerEdge R720
machine with two Intel Xeon E5-2690 CPUs, 16 cores and 32 threads in total,
2.9-3.8~GHz, and 192~GB of RAM running the Gentoo Linux operating system.

\begin{table*}
\centering
\begin{tabular}{|c|c|c||r|r|r|r|r|r|}
\hline
\multicolumn{3}{|c||}{Example} & \multicolumn{1}{c|}{Magma} &
  \multicolumn{5}{c|}{Singular} \\
\hline
\multirow{2}{*}{Ideal} & Min. & deg. &
  \multicolumn{1}{c|}{\texttt{Groebner}} &
  \multicolumn{1}{c|}{\multirow{2}{*}{\texttt{std}}} &
  \multicolumn{2}{c|}{\texttt{modStd}} &
  \multicolumn{2}{c|}{\texttt{nfmodStd}} \\
\cline{6-9}
& Poly. & of $m_i$ & \multicolumn{1}{c|}{\texttt{Basis}} & & 1 c. & 32 c. &
  1 c. & 32 c. \\
\hline\hline
I1  & $m_1$ &  2 & 1241.98 &     1.51 &     1.24 &    0.37 &   0.22 &  0.13 \\
\hline
I2  & $m_2$ &  5 &   error &    70.55 &    19.59 &    4.79 &   1.89 &  0.61 \\
\hline
I3a & $m_3$ &  7 &       - &     0.90 &   143.79 &    9.34 &   3.27 &  0.51 \\
\hline
I3b & $m_3$ &  7 &       - &   314.00 & 11212.00 & 1118.78 &  97.43 & 19.23 \\
\hline
I4  & $m_4$ &  6 &       - &   265.53 &  9163.38 &  567.03 & 686.01 & 99.41 \\
\hline
I5  & $m_5$ & 12 &       - &  2061.95 &  3321.28 &  256.58 & 430.23 & 71.47 \\
\hline
I6  & $m_6$ &  2 &    2.93 &  8931.13 &   197.20 &   47.54 &  24.26 &  8.99 \\
\hline
I7  & $m_7$ &  8 &       - &     0.90 &  2044.08 &  195.41 &   8.54 &  1.87 \\
\hline
I8  & $m_8$ &  7 &       - & 15477.87 & 15274.97 & 4787.49 &  92.99 & 23.89 \\
\hline
\end{tabular}
\captionsetup{width=.9\textwidth}
\caption{Total running times in seconds for computing a Gr\"obner basis of the
considered ideals with the corresponding minimal polynomial via
\textnormal{\texttt{GroebnerBasis}}, \textnormal{\texttt{std}},
\textnormal{\texttt{modStd}} and \textnormal{\texttt{nfmodStd}}, using 1~core
and 32~cores where applicable}%
\label{tab:timings}
\end{table*}

The results are summarized in Table~\ref{tab:timings}. Some of the computations
in Magma did not finish within 12~hours. This is indicated by a dash (-). Note
that in all those cases, the computation also occupied an excessive amount of
memory, more than 100~GB at the point when we interrupted it. All timings are
in seconds. We use the degree reverse lexicographical ordering (\verb+dp+ in
\textsc{Singular}) for all examples.

In our implementation, the number of primes which are chosen in
line~\ref{line:modStd_variant_primes} of Algorithm~\ref{alg:modStd_variant}
depends on the number of cores. For our timings, we started with 10 primes on
one core and 25 primes on 32 cores. The runtime depends heavily on the
splitting behaviour of the minimal polynomial modulo the chosen primes. Finding
the optimal strategy for this is still under active research.

\begin{rem}
We understand that Magma has no parallel version of the Gr\"obner basis
algorithm which works over algebraic number fields. Therefore we have conducted
the timings in Magma using one core only.
\end{rem}

From Table~\ref{tab:timings}, we see that the \textsc{Singular} commands
\texttt{std} and \texttt{modStd} perform well in comparison to the Magma
command \texttt{GroebnerBasis}. However, one can see that our algorithm
\texttt{nfmodStd} is even much faster.

\section{Acknowledgments}

We would like to thank Gerhard Pfister for many fruitful discussions.

\bibliographystyle{abbrv}
\bibliography{nfmodstd}

\section*{Appendix}

The following are the benchmark problems used to demonstrate the efficiency of
Algorithm~\ref{alg:nfmodStd}. They are available in the source code of the
\textsc{Singular} library \texttt{nfmodstd.lib} \cite{nfmodstd}.

\begin{enumerate}
\addtolength{\itemsep}{0pt plus 2pt minus 2pt}
\item
\begin{verbatim}
ring R = (0,a), (x,y,z), dp;
minpoly = (a^2+1);
poly f1 = (a+8)*x^2*y^2+5*x*y^3+(-a+3)*x^3*z+x^2*y*z;
poly f2 = x^5+2*y^3*z^2+13*y^2*z^3+5*y*z^4;
poly f3 = 8*x^3+(a+12)*y^3+x*z^2+3;
poly f4 = (-a+7)*x^2*y^4+y^3*z^3+18*y^3*z^2;
ideal I1 = f1,f2,f3,f4;
\end{verbatim}
\item
\begin{verbatim}
ring R = (0,a), (x,y,z), dp;
minpoly = (a^5+a^2+2);
poly f1 = 2*x*y^4*z^2+(a-1)*x^2*y^3*z+(2*a)*x*y*z^2+7*y^3
          +(7*a+1);
poly f2 = 2*x^2*y^4*z+(a)*x^2*y*z^2-x*y^2*z^2+(2*a+3)*x^2*y*z
          -12*x+(12*a)*y;
poly f3 = (2*a)*y^5*z+x^2*y^2*z-x*y^3*z+(-a)*x*y^3+y^4
          +2*y^2*z;
poly f4 = (3*a)*x*y^4*z^3+(a+1)*x^2*y^2*z-x*y^3*z+4*y^3*z^2
          +(3*a)*x*y*z^3+4*z^2-x+(a)*y;
ideal I2 = f1,f2,f3,f4;
\end{verbatim}
\item
\begin{verbatim}
ring R = (0,a), (v,w,x,y,z), dp;
minpoly = (a^7-7*a+3);
poly f1 = (a)*v+(a-1)*w+x+(a+2)*y+z;
poly f2 = v*w+(a-1)*w*x+(a+2)*v*y+x*y+(a)*y*z;
poly f3 = (a)*v*w*x+(a+5)*w*x*y+(a)*v*w*z+(a+2)*v*y*z
          +(a)*x*y*z;
poly f4 = (a-11)*v*w*x*y+(a+5)*v*w*x*z+(a)*v*w*y*z+(a)*v*x*y*z
          +(a)*w*x*y*z;
poly f5 = (a+3)*v*w*x*y*z+(a+23);
ideal I3a = f1,f2,f3,f4,f5;
\end{verbatim}
\item
\begin{verbatim}
ring R = (0,a), (u,v,w,x,y,z), dp;
minpoly = (a^7-7*a+3);
poly f1 = (a)*u+(a+2)*v+w+x+y+z;
poly f2 = u*v+v*w+w*x+x*y+(a+3)*u*z+y*z;
poly f3 = u*v*w+v*w*x+(a+1)*w*x*y+u*v*z+u*y*z+x*y*z;
poly f4 = (a-1)*u*v*w*x+v*w*x*y+u*v*w*z+u*v*y*z+u*x*y*z
          +w*x*y*z;
poly f5 = u*v*w*x*y+(a+1)*u*v*w*x*z+u*v*w*y*z+u*v*x*y*z
          +u*w*x*y*z+v*w*x*y*z;
poly f6 = u*v*w*x*y*z+(-a+2);
ideal I3b = f1,f2,f3,f4,f5,f6;
\end{verbatim}
\item
\begin{verbatim}
ring R = (0,a), (w,x,y,z), dp;
minpoly = (a^6+a^5+a^4+a^3+a^2+a+1);
poly f1 = (a+5)*w^3*x^2*y+(a-3)*w^2*x^3*y+(a+7)*w*x^2*y^2;
poly f2 = (a)*w^5+(a+3)*w*x^2*y^2+(a^2+11)*x^2*y^2*z;
poly f3 = (a+7)*w^3+12*x^3+4*w*x*y+(a)*z^3;
poly f4 = 3*w^3+(a-4)*x^3+x*y^2;
ideal I4 = f1,f2,f3,f4;
\end{verbatim}
\item
\begin{verbatim}
ring R = (0,a), (w,x,y,z), dp;
minpoly = (a^12-5*a^11+24*a^10-115*a^9+551*a^8-2640*a^7
          +12649*a^6-2640*a^5+551*a^4-115*a^3+24*a^2-5*a+1);
poly f1 = (2*a+3)*w*x^4*y^2+(a+1)*w^2*x^3*y*z+2*w*x*y^2*z^3
          +(7*a-1)*x^3*z^4;
poly f2 = 2*w^2*x^4*y+w^2*x*y^2*z^2+(-a)*w*x^2*y^2*z^2
          +(a+11)*w^2*x*y*z^3-12*w*z^6+12*x*z^6;
poly f3 = 2*x^5*y+w^2*x^2*y*z-w*x^3*y*z-w*x^3*z^2+(a)*x^4*z^2
          +2*x^2*y*z^3;
poly f4 = 3*w*x^4*y^3+w^2*x^2*y*z^3-w*x^3*y*z^3
          +(a+4)*x^3*y^2*z^3+3*w*x*y^3*z^3+(4*a)*y^2*z^6-w*z^7
          +x*z^7;
ideal I5 = f1,f2,f3,f4;
\end{verbatim}
\item
\begin{verbatim}
ring R = (0,a), (u,v,w,x,y,z), dp;
minpoly = (a^2+5*a+1);
poly f1 = u+v+w+x+y+z+(a);
poly f2 = u*v+v*w+w*x+x*y+y*z+(a)*u+(a)*z;
poly f3 = u*v*w+v*w*x+w*x*y+x*y*z+(a)*u*v+(a)*u*z+(a)*y*z;
poly f4 = u*v*w*x+v*w*x*y+w*x*y*z+(a)*u*v*w+(a)*u*v*z
          +(a)*u*y*z+(a)*x*y*z;
poly f5 = u*v*w*x*y+v*w*x*y*z+(a)*u*v*w*x+(a)*u*v*w*z
          +(a)*u*v*y*z+(a)*u*x*y*z+(a)*w*x*y*z;
poly f6 = u*v*w*x*y*z+(a)*u*v*w*x*y+(a)*u*v*w*x*z
          +(a)*u*v*w*y*z+(a)*u*v*x*y*z+(a)*u*w*x*y*z
          +(a)*v*w*x*y*z;
poly f7 = (a)*u*v*w*x*y*z-1;
ideal I6 = f1,f2,f3,f4,f5,f6,f7;
\end{verbatim}
\item
\begin{verbatim}
ring R = (0,a), (w,x,y,z), dp;
minpoly = (a^8-16*a^7+19*a^6-a^5-5*a^4+13*a^3-9*a^2+13*a+17);
poly f1 = (-a^2-1)*x^2*y+2*w*x*z-2*w+(a^2+1)*y;
poly f2 = (a^3-a-3)*w^3*y+4*w*x^2*y+4*w^2*x*z+2*x^3*z+(a)*w^2
          -10*x^2+4*w*y-10*x*z+(2*a^2+a);
poly f3 = (a^2+a+11)*x*y*z+w*z^2-w-2*y;
poly f4 = -w*y^3+4*x*y^2*z+4*w*y*z^2+2*x*z^3+(2*a^3+a^2)*w*y
          +4*y^2-10*x*z-10*z^2+(3*a^2+5);
ideal I7 = f1,f2,f3,f4;
\end{verbatim}
\item
\begin{verbatim}
ring R = (0,a), (t,u,v,w,x,y,z), dp;
minpoly = (a^7+10*a^5+5*a^3+10*a+1);
poly f1 = v*x+w*y-x*z-w-y;
poly f2 = v*w-u*x+x*y-w*z+v+x+z;
poly f3 = t*w-w^2+x^2-t;
poly f4 = (-a)*v^2-u*y+y^2-v*z-z^2+u;
poly f5 = t*v+v*w+(-a^2-a-5)*x*y-t*z+w*z+v+x+z+(a+1);
poly f6 = t*u+u*w+(-a-11)*v*x-t*y+w*y-x*z-t-u+w+y;
poly f7 = w^2*y^3-w*x*y^3+x^2*y^3+w^2*y^2*z-w*x*y^2*z
          +x^2*y^2*z+w^2*y*z^2-w*x*y*z^2+x^2*y*z^2+w^2*z^3
          -w*x*z^3+x^2*z^3;
poly f8 = t^2*u^3+t^2*u^2*v+t^2*u*v^2+t^2*v^3-t*u^3*x
          -t*u^2*v*x-t*u*v^2*x-t*v^3*x+u^3*x^2+u^2*v*x^2
          +u*v^2*x^2+v^3*x^2;
ideal I8 = f1,f2,f3,f4,f5,f6,f7,f8;
\end{verbatim}
\end{enumerate}

\clearpage

\end{document}